\documentclass[amstex]{article}
\usepackage{amsmath,amsfonts}
\input{tcilatex}

\begin{document}

\nocite{*}

\title{Differential and complex geometry of two-dimensional noncommutative tori}

\author{Momar Dieng\\ Albert Schwarz\thanks{Partially supported by NSF grant No. DMS 9971304}\\ \\
Department of Mathematics\\ University of California \\ Davis, CA 95616}

\date{}

\maketitle

\begin{abstract}
We analyze in detail projective modules over two-dimensional noncommutative
tori \ and complex structures on these modules.We concentrate our attention
on properties of holomorphic vectors in these modules; the theory of these
vectors generalizes the theory of theta-functions. The paper is
self-contained; it can be used also as an introduction to the theory of
noncommutative spaces with simplest space of this kind thoroughly analyzed
as a basic example.
\end{abstract}

\section{Introduction}

Differential geometry of noncommutative spaces, in particular, of
noncommutative tori was developed by A. Connes \cite{Co1}, \cite{Co2}. More detailed
analysis of differential geometry of noncommutative tori was performed in
%\cite{Co3}--\cite{Ko2}
\cite{Co3}, \cite{Co4},\cite{Ri1}, \cite{Ri2}, \cite{Sc1},
\cite{Sc2},\cite{Ko2}.

Partly, the interest to this subject was motivated by the applications to
string/M-theory found in \cite{Co4} (see \cite{Ko3} for review). Complex geometry of noncommutative tori and of
projective modules over them was studied in \cite{Sc3} in connection with
noncommutative generalization of theta-functions.

The goal of present paper is to illustrate general results of \cite{Sc3} using the
example of two-dimensional noncommutative tori, where all calculations can
be performed explicitly. We repeat all basic definitions, starting with the
definition of noncommutative torus. Therefore, the paper can be read
independently of \cite{Sc3} and of other papers about noncommutative tori. However,
in the rest of the introduction, we assume some knowledge of preceding work.

In Sec. 2 we describe projective modules over two-dimensional
noncommutative
torus and their tensor products. The results of this section are closely
related to some results proved in multidimensional case in \cite{Sc3}. We will
formulate these general results restricting ourselves to basic $T_{\theta }$%
- modules (modules with constant curvature connection where the endomorphism
algebra is another noncommutative torus $T_{\theta ^{\prime }}$). Instead of
basic modules we can work with corresponding $(T_{\theta ^{\prime
}},T_{\theta })$-bimodules. Tensor product of $(T_{\theta ^{\prime
}},T_{\theta })$-bimodule ${\cal E}_{1}$ and $(T_{\theta },T_{\theta
^{\prime \prime }})$-bimodule ${\cal E}_{2}$ is a $(T_{\theta ^{\prime
}},T_{\theta ^{\prime \prime }})$-bimodule ${\cal E}.$ We can find the
bimodule ${\bf \ }{\cal E}$\ using the fact that basic bimodules ${\cal E}%
_{1}$ and ${\cal E}_{2}$ can be described by means of elements $g_{1},g_{2}\in SO(d,d,{\mathbb Z)}$ obeying $\theta ^{\prime }=g_{1}\theta
,\theta ^{\prime \prime }=g_{2}\theta$. Here $g_{1},g_{2}$ act on the
parameter of noncommutativity - antisymmetric matrix $\theta $ - by means of
fractional linear transformation. The description of ${\cal \ \ E}$
follows from the relation $\theta ^{\prime \prime }=g_{2}g_{1}^{-1}\theta
^{\prime }$. In two-dimensional case all basic right $T_{\theta }$--modules can be represented as modules $E_{n,m}$ (see Sec. 2) where $m$ and
$n$ are relatively prime. The space $E_{n,m}$ can be considered also as a left $%
T_{\theta ^{\prime }}$-module $E_{a,m}^{\prime }$ where $%
an-bm=1,a,b\in {\mathbb Z,\theta }^{\prime }=(a\theta +b)(m\theta +n)^{-1}$, or
as $(T_{\theta ^{\prime }},T_{\theta })$-bimodule that will be denoted by $%
{\cal E}_{1}$. One can use general results to calculate the tensor product $%
{\cal E}_{1}\underset{T_{\theta }}{\otimes }{\cal E}_{2}$ where ${\cal E}%
_{2} $ is a $(T_{\theta },T_{\theta ^{\prime \prime }})$-bimodule,
corresponding to the basic left $T_{\theta }$-module $E_{k,l}^{\prime }$. It follows from these results that this tensor product considered as a left $%
T_{\theta ^{\prime }}$-module is isomorphic to
$E_{ak+bl,nl+mk}^{\prime }$. We give a direct proof of this
relation constructing an explicit isomorphism. This construction
is used essentially in Sec. 3, devoted to complex geometry of
projective modules over two-dimensional noncommutative tori. We
calculate holomorphic vectors (theta-vectors) in basic modules \
and tensor products of these vectors. The relations we obtain
generalize well known relations for theta-functions. The appendix
contains a very detailed calculation of the tensor product (we did
not include the appendix into the journal version). We did not try
to analyze the connection of our results to Manin's version of the
theory of noncommutative theta--functions \cite{Ma1}--\cite{Ma3}.
This is an interesting problem.

\section{\protect\bigskip Projective modules and their tensor products.}

One can define the algebra $T_{\theta }^{d}$ of smooth functions on $d$%
-dimensional noncommutative torus as an algebra of formal linear
combinations $f=\sum C_{\overrightarrow{n}}U_{\overrightarrow{n}}$ where $%
\overrightarrow{n}\in \mathbb Z^{d}$ and $C_{\overrightarrow{n}}$ are complex
numbers tending to zero at infinity faster then any power and the
multiplication is governed by the rule

$$U_{\overrightarrow{n}}U_{\overrightarrow{m}}=e^{\pi i\overrightarrow{n}
\theta \overrightarrow{m}}U_{\overrightarrow{n}+\overrightarrow{m}} $$

\noindent where $\theta =\theta ^{\alpha \beta }$ is a real antisymmetric
matrix. An antilinear involution on $T_{\theta }^{d}$ is defined by the
requirement $U_{\overrightarrow{n}}^{\ast }=U_{-\overrightarrow{n}}$. These
operations together with the standard structure of vector space permit us to
consider $T_{\theta }^{d}$ as involutive associative algebra with unit
element $1=U_{\vec{0}}$ and with trace Tr$f=C_{\overrightarrow{0}%
} $. We will fix our attention on the case $d=2;$ in this case the matrix $%
\theta ^{\alpha \beta }$ can be specified by means of one number $\theta
^{12}=\theta $ and the multiplication is specified by the relation

$$U_{1}U_{2}=e^{2\pi i\theta }U_{2}U_{1} $$

\noindent where $U_{\alpha }=U_{\overrightarrow{e_{\alpha }}}$ are elements $U_{%
\overrightarrow{n}}$ corresponding to vectors of standard basis $%
\overrightarrow{e}_{1}=(1,0)$, $\overrightarrow{e}_{2}=(0,1)$.

In the case when $\theta $ is a natural number $T_{\theta }^{2}$ is
isomorphic to the algebra of smooth functions on the commutative
two-dimensional torus. We will consider projective modules over $T_{\theta
}^{2}$ (we always assume that our modules are finitely generated; a
projective module is by definition a direct summand in a free module $%
(T_{\theta }^{2})^{k}$). We assume that the number $\theta $ is irrational;
in this case every projective right module over $T_{\theta }^{2}$ is
isomorphic to one of the modules $E_{n,m}$ defined in the following way (see
\cite{Co2}, \cite{Co3}, \cite{Ko3}). The elements of $E_{n,m}$ are functions on the Schwartz space ${\cal S}({\mathbb R}\times{\mathbb Z}_{m})$ and the action of the generators $U_{1}, U_{2}$ of $T_{\theta }^{2}$ is given by:

$$U_{1}\,f(x, \mu)=f(x-\frac{n+m\,\theta}{m}, \mu-1)$$
$$ U_{2}\,f(x, \mu)=e^{2\,\pi\,i(x-\mu\,n/m)}\,f(x, \mu)$$

\noindent We can define also left modules over $T_{\theta }^{2}$ replacing \ $\theta
$ by $-\theta $ in the above formulas; we will use the notation $E_{n,m}^{\prime }$ for left modules. If it is necessary to emphasize that we consider $E_{n, m}$ or $
E_{n,m}^{\prime }$ as a $T_{\theta }$--modules we will use the notation $E_{n, m}(\theta)$ or $E_{n, m}^{\prime}(\theta)$.
We will assume that the numbers $m$ and $n$ are relatively prime.
Corresponding modules $E_{n,m}$ are called basic modules. Every module over $T^{2}_{\theta}$ can be represented as a direct sum of isomorphic basic modules. In this case the algebra of
endomorphisms of the module $E=E_{n,m}$ is again a noncommutative torus $%
T_{\theta ^{\prime }}$ generated by operators
$$Z_{1}\,f(x, \mu)=f(x-\frac{1}{m}, \mu-a)$$
$$ Z_{2}\,f(x, \mu)=\exp\left[2\,\pi\,i\left(\frac{x}{n+m\theta}-\frac{\mu}{m}\right)\right]\,f(x, \mu)$$

\noindent where $a,b$ are integers obeying $an-bm=1$. This fact follows from relation
$$Z_{1}Z_{2}=e^{-2\pi i\theta ^{\prime }}Z_{2}Z_{1} $$

\noindent where $\theta ^{\prime
}=\frac{b+a\,\theta}{n+m\,\theta}.$ One can consider $E$ as a left
$T_{\theta ^{\prime }}$-module; this module is isomorphic to
$E_{a,m}^{\prime }.$

We can regard $E$ also as $(T_{\theta ^{\prime }},T_{\theta })$-bimodule, since the action of $T_{\theta ^{\prime }}$ commutes with the action of $T_{\theta}$. Notice that the bimodule $E$ depends on the choice of $a$ and $b$. The tori $T_{\theta }$ and $T_{\theta ^{\prime }}$ are Morita equivalent and the bimodule $E$ is a Morita equivalence bimodule (\cite{Co2}, \cite{Co3}, \cite{Ko3}).

Recall that two associative algebras $A$ and $B$ are Morita equivalent if
corresponding categories of modules are equivalent. Having an $(A,B)$%
-bimodule $P$ we can assign to every right $A$-module $E$ a right
$B$-module $E^{\prime }=E\underset{A}{\otimes }P$ where tensor
product over $A$ is obtained  from the standard tensor product
over ${\mathbb C}$ by means of identification $ea\otimes p\sim
e\otimes ap$. If this correspondence is invertible we say that $P$
is a Morita equivalence bimodule.

 Let us calculate now the tensor product $E_{n,m}\underset{%
T_{\theta }}{\otimes }E_{k,l}^{\prime }$ where $E_{n,m}$ is a right $%
T_{\theta }^{2}$-module and $E_{k,l}^{\prime }$ is a left $T_{\theta }^{2}$%
-module. We will see that this product can be considered a vector space $%
{\cal E=S}({\mathbb R}\times {\mathbb Z}_{nl+mk})$. Namely, the formula:
\begin{eqnarray}h(z,\Delta)= \sum_{q\in\mathbb Z}\left[f\left((n+m\,\theta)z-\frac{n+m\,\theta}{m}q+\frac{l(n+m\,\theta)}{m(n\,l+m\,k)}\Delta, -q+a\,\Delta\right)\right.\nonumber\\ \left.\cdot g\left((n+m\,\theta)z+\frac{k-l\,\theta}{l}q-\frac{k-l\,\theta}{n\,l+m\,k}\Delta, q\right) \right]\label{map}
\end{eqnarray}

\noindent where $0\leq \Delta\leq nl+mk-1$, specifies a map of tensor product over $E_{n,m}\underset{T_{\theta}}{\otimes }E_{k,l}^{\prime }$ onto ${\cal E}$, so that, with the usual notation, $f\in E_{n,m}(\theta)$ , $g\in E^{\prime}_{k,l}(\theta)$ and $h \in E^{\prime}_{a\,k+b\,l,n\,l+m\,k}(\theta^{\prime})$. This map is an isomorphism
between the tensor product and ${\cal E}$. The formula (\ref{map}) determines a
bilinear map that is compatible with identification $ea\otimes p\sim
e\otimes ap$ (see Appendix for details). We can consider $E_{n,m}$ as a $%
(T_{\theta ^{\prime }},T_{\theta })$-bimodule where $\theta ^{\prime }=%
\frac{b+a\theta }{n+m\theta },an-bm=1$ and $E_{k,l}$ as a $(T_{\theta
},T_{\theta ^{\prime \prime }})$-bimodule where $\theta ^{\prime \prime }=-%
\frac{d-c\theta }{k-l\theta }$, $ck-dl=1$. Then ${\cal E}$ can be regarded
as $(T_{\theta ^{\prime }},T_{\theta ^{\prime \prime }})$-bimodule. It
follows from (\ref{map}) that ${\cal E}$ considered as a left $T_{\theta ^{\prime }}$%
-module is isomorphic to $E_{N^{\prime },M}^{\prime }$ where $%
M=nl+mk,N^{\prime }=ak+bl$ and ${\cal E}$ considered as a right $T_{\theta
^{\prime \prime }}$-module is isomorphic to $E_{N^{\prime \prime },M}$ where
$N^{\prime \prime }=-(cn+md)$ for $ck-dl=1$. This result can be obtained also from
general considerations (see Introduction).

Let us apply (\ref{map}) to the special case when:

$$f(x,\mu)=e^{-\frac{1}{2}\sigma_{_{1}} x^{2}-c_{_{1}}x}\delta_{\alpha}^{\mu }$$
$$g(y,\nu)=e^{-\frac{1}{2}\sigma_{_{2}} y^{2}-c_{_{2}}x}\delta_{\beta}^{\nu }$$

\noindent where
$$\sigma_{_{1}} =i\tau_{1} \frac{m}{n+m\,\theta },\qquad \mu\in {\mathbb Z}_{m}$$
$$\sigma_{_{2}} =i\tau_{2} \frac{l}{k-l\,\theta },\qquad \nu\in {\mathbb Z}_{l}$$

\noindent and $\delta_{j}^{i}$ is the usual Kronecker--delta. It is straightforward to check that in this case:
\begin{align}
h_{\alpha\,\beta}(z,\Delta)= \sum_{q\in\mathbb Z}\left[\exp\left\{-\frac{\sigma_{_{1}}}{2}\left((n+m\,\theta)z-\frac{n+m\,\theta}{m}q+\frac{l(n+m\,\theta)}{m(n\,l+m\,k)}\Delta\right)^{2}\right.\right.\nonumber\\
\left.-c_{_{1}}\left((n+m\,\theta)z-\frac{n+m\,\theta}{m}q+\frac{l(n+m\,\theta)}{m(n\,l+m\,k)}\Delta\right)\right\}\,\delta^{-q-a\,\Delta}_{\alpha}\nonumber\\
 \left.\cdot \exp\left\{-\frac{\sigma_{_{2}}}{2}\left((n+m\,\theta)z+\frac{k-l\,\theta}{l}q-\frac{k-l\,\theta}{n\,l+m\,k}\Delta\right)^{2}\right.\right.\nonumber\\
\left.\left.-c_{_{2}}\left((n+m\,\theta)z+\frac{k-l\,\theta}{l}q-\frac{k-l\,\theta}{n\,l+m\,k}\Delta\right)\right\}\,\delta^{q}_{\beta} \right]\label{prodbasis}
\end{align}

\noindent Notice that the right hand side of (\ref{prodbasis}) can be expressed in terms of theta--functions; first note that the set of solutions to the system of congruences:
$$\begin{cases}
q=a\,\Delta-\alpha \pmod m\\
q= \beta \pmod l
\end{cases}$$

\noindent can be written as $q_{o}+u\frac{m\,l}{r}$ for some integer $q_{o}$, $r=\gcd(m,l)$, and $u\in\mathbb Z$. Hence we can substitute $q_{o}+u\frac{m\,l}{r}$ in for $q$ in (\ref{prodbasis}), do away with Kronecker deltas, and sum over $u$ instead of $q$. We obtain:

$$h_{\alpha\,\beta}(z,\Delta)=\Theta(s,t)\cdot\xi(z,\Delta)$$

\noindent where

\begin{eqnarray*}
\xi_{\alpha\,\beta}(z,\Delta)&=&\exp\left[ -\frac{\sigma_{_{1}}}{2}\left((n+m\,\theta)\,z+\frac{l(n+m\,\theta)\,\Delta}{m(n\,l+m\,k)}\right)^{2}-\frac{\sigma_{_{2}}}{2}\left((n+m\,\theta)\,z-\frac{(k-l\,\theta)\,\Delta}{n\,l+m\,k}\right)^{2}\right.\\
&& \qquad \qquad -(c_{_{1}}+c_{_{2}})\,(n+m\,\theta)\,z-\left(c_{_{1}}\frac{l(n+m\,\theta)}{m(n\,l+m\,k)}-c_{_{2}}\frac{k-l\,\theta}{n\,l+m\,k}\right)\,\Delta \\
&& \qquad\qquad \qquad \left. -\left(\sigma_{_{1}}\left(\frac{n+m\,\theta}{m}\right)^{2}+\sigma_{_{2}}\left(\frac{k-l\,\theta}{l}\right)^{2}\right)\frac{q_{o}^{2}}{2} \right]
\end{eqnarray*}

%\begin{eqnarray*}
%\xi_{\alpha\,\beta}(z,\Delta)&=&\exp\left[\frac{\pi\,i\,\tau(n\,l+m\,k)(n+m\,\theta)}{k-l\,\theta}z^{2}-(c_{_{1}}+c_{_{2}})(n+m\,\theta)z\right.\\ && \qquad\qquad \left.-\frac{\pi\,i\,\tau\,l}{m(n\,l+m\,k)}\,\Delta^{2}-\frac{\left(l(n+m\,\theta)c_{_{1}}-m(k-l\,\theta)c_{_{2}}\right)}{m(n\,l+m\,k)}\,\Delta\right]
%\end{eqnarray*}

$$s=-\frac{\sigma_{_{1}}\,l^{2}\,(n+m\,\theta)^{2}+\sigma_{_{2}}\,m^{2}\,(k-l\,\theta)^{2}}{2\,\pi\,i\,r^{2}}$$

%$$s=-m\,l\,(n\,l+m\,k)$$

\begin{eqnarray*}t &=&\frac{\sigma_{_{1}}\,(n+m\,\theta)^{2}}{2\,\pi\,i\,m}\left(z+\frac{l\,\Delta}{m\,(n\,l+m\,k)}-\frac{l\,q_{o}}{r} \right) + c_{_{1}}\,\frac{n+m\,\theta}{2\,\pi\,i\,m}- c_{_{2}}\,\frac{k-l\,\theta}{2\,\pi\,i\,l} \\
&&\qquad +\frac{\sigma_{_{2}}\,(k-l\,\theta)^{2}}{2\,\pi\,i\,l}\left(-\frac{n+m\,\theta}{k-l\,\theta}\,z+\frac{\Delta}{n\,l+m\,k}-\frac{m\,q_{o}}{r} \right)
\end{eqnarray*}

%$$t=\frac{l\,\tau\,\Delta-q_{o}(n\,l+m\,k)}{r} + \frac{c_{_{1}}\,l\,(n+m\,\theta)-c_{_{2}}\,m\,(k-l\,\theta)}{2\pi\,i\,r}$$

\noindent and $\Theta(s,t)$ a classical theta--function. We recall that classical--theta functions are defined as:

$$\Theta(s,t)=\sum_{u\in \mathbb Z}e^{\pi\,i\,s\,u^{2}+2\pi\,i\,t\,u}\quad , \qquad \textrm{Im }s>0$$

\section{ Connections and complex structures.}

One can define a connection on a right $T_{\theta }^{d}$-module
$E$ as a set of ${\mathbb C}$-linear operators $\nabla
_{1},...,\nabla _{d}$ obeying Leibniz rule:
\begin{equation}
\nabla _{\alpha }(ea)=\nabla _{\alpha }e\cdot a+e\cdot \delta _{\alpha }a\label{con}
\end{equation}

\noindent where $e\in E,a\in T_{\theta }^{d}$ and the derivatives $\delta _{\alpha }a$
are specified by the formula $\delta _{\alpha }U_{\overrightarrow{n}}=2\pi
niU_{\overrightarrow{n}}$ (see[1]). It will be convenient for us to
generalize the notion of connection replacing $\delta _{\alpha }$ in (\ref{con}) by $%
\delta _{\alpha }^{\prime }=a_{\alpha }^{\beta }\delta _{\beta }$ where $%
a_{\alpha }^{\beta }$ is a non-degenerate matrix. A connection $\nabla
_{1},...,\nabla _{d}$ is a constant curvature connection if
$$\left[ \nabla _{\alpha },\nabla _{\beta }\right]=i\,f_{\alpha \beta }\cdot 1$$
\noindent where $f_{\alpha \beta }$ are numbers and $1$ stands for
identity operator. Similar definitions can be given for left
modules. We always consider unitary connections (i.e. operators
$\nabla _{\alpha }$ should be anti\ Hermitian). It is easy to
check that operators
\begin{equation}
\nabla _{1}=2\,\pi\,i\frac{m}{n+m\theta }x,\text{ \ \ }\nabla
_{2}=2\,\pi\,\frac{d}{dx}\label{conscon}
\end{equation}

\noindent specify a constant curvature connection of right $T_{\theta }$-module $%
E_{m,n}.$ The same operators determine a connection (in
generalized sense) on $E_{n,m}$ considered as a left $T_{\theta
^{\prime }}$-module. This follows from relations $[\nabla
_{i},\nabla_{j}]=\frac{2\pi\,i\,m}{n+m\theta}$. We see that
$\nabla _{1}, \nabla _{2}$ can be considered as a constant
curvature connection on $(T_{\theta ^{\prime }},T_{\theta
})$-bimodule.

If $\nabla _{1},\nabla _{2}$ is a constant curvature connection on
a module $E$ we can introduce a complex structure on $E$ fixing
$\overline{\partial \text{-}}$connection \ $\overline{\nabla
}=\lambda_{1} \nabla _{1}+\lambda_{2} \nabla _{2},$where
$\lambda_{1} $ and $\lambda_{2}$ are complex numbers and the
quotient $\tau =\lambda_{1}/\lambda_{2} $ is not real \cite{Sc3}.
This complex structure on a $T_{\theta } $-module corresponds to
complex structure on $T_{\theta },$ specified by means of
$\overline{\partial }$-derivative $\lambda_{1} \delta
_{1}+\lambda_{2} \delta _{2}$. By definition, vector $\Theta \in
E$ is holomorphic if:
$$\overline{\nabla }\Theta =0$$

\noindent Notice that the notion of holomorphicity depends only on
$\tau =\lambda_{1}/\lambda_{2} $, therefore we say that
$\overline{\nabla }$ and $\rho \overline{\nabla }$ where $\rho
\neq 0$ determine the same complex structure on $E.$\ Holomorphic
vectors are closely related to theta-functions, hence we call \
holomorphic vectors in basic modules theta-vectors.

Let us consider \ holomorphic vectors in $T_{\theta }$-modules $E_{n,m},$
assuming that $n$ and $m$ are relatively prime. In this case all\ constant
curvature connections have the form $\nabla _{\alpha }=\nabla _{\alpha
}^{0}+c_{\alpha }$, where $\nabla _{\alpha }^{0}$ stand for the connection
(\ref{conscon}) and $c_{_{1}},c_{_{2}}$ are constants. The equation $\overline{\nabla }\Theta =0
$ takes the form
$$(i\tau \frac{m}{n+m\theta }x+\frac{\partial }{\partial x}+c)\varphi (x,\mu
)=0$$

\noindent If $\limfunc{Im}\tau <0$ it has $m$ linearly independent solutions
\begin{equation}
\varphi _{\alpha }(x,\mu )=e^{-\frac{1}{2}\sigma x^{2}-cx}\delta _{\alpha
}^{\mu }\label{stdbasis}
\end{equation}

\noindent where
$$\sigma =i\tau \frac{m}{n+m\,\theta },\qquad \mu\in {\mathbb Z}_{m}$$

\noindent The functions (\ref{stdbasis}) belong to ${\cal S}$ only if $\limfunc{Re}\sigma >0$. We
assumed that $n+m\theta >0$; in the case $n+m\theta <0$, the condition of
existence of holomorphic vectors is that $\limfunc{Im}\tau >0$. We see that in
the case when holomorphic vectors exist the space ${\cal H}_{n,m}$ of holomorphic vectors in $E_{n,m}$ is $m$-dimensional; the functions (\ref{stdbasis})
constitute a basis of ${\cal H}_{n,m}$. Considering $E_{n,m}$ as a $(T_{\theta ^{\prime }},T_{\theta })$-bimodule and taking into account that a
constant curvature connection on $E_{n,m}$ is a constant curvature connection on the
bimodule, one can define a notion of complex structure and of holomorphic vector
for a bimodule. More precisely, complex structure on $T_{\theta }$-module $E_{n,m}$ induces a complex structure on same space considered as left $T_{\theta ^{\prime }}$-module in such a way that the notion of holomorphic vector remains the same. These two complex structures specify a
complex structure on a bimodule; the notion of holomorphic vector in a
bimodule coincides with corresponding notion for both modules. If ${\cal E}^{\prime }$ is a complex  $(T_{\theta ^{\prime}},T_{\theta })$-bimodule, ${\cal E}^{\prime \prime }$ is a complex  $(T_{\theta },T_{\theta ^{\prime \prime }})$-bimodule, then the  $(T_{\theta
^{\prime }},T_{\theta ^{\prime \prime }})$-bimodule ${\cal E}={\cal E}^{\prime }\underset{T_{\theta }}{\otimes }{\cal E}^{\prime \prime }$ can be
equipped with complex structure. We assume that  complex structure on the
right $T_{\theta }$-module ${\cal E}^{\prime }$ and on the left $T_{\theta }$%
-module ${\cal E}^{\prime \prime }$ correspond to the same complex structure
on $T_{\theta }$. One can prove \cite{Sc3} that tensor product of two \
holomorphic vectors is again a  holomorphic vector (i.e. we have a natural
map ${\cal H}^{\prime }\otimes {\cal H}^{\prime \prime }$ into ${\cal H}$
where ${\cal H}^{\prime },{\cal H}^{\prime \prime }$ and ${\cal H}$ stand
for spaces of holomorphic vectors in ${\cal E}^{\prime },{\cal E}^{\prime
\prime }$ and ${\cal E}$ correspondingly).
Let the basis of ${\cal H}^{\prime }$ consist of:

$$\varphi^{\prime}_{\alpha}=e^{-\frac{1}{2}\sigma_{_{1}} x^{2}-c_{_{1}}x}\delta_{\alpha}^{\mu }$$

\noindent where

$$\sigma_{_{1}} =i\tau_{1} \frac{m}{n+m\,\theta },\qquad \mu\in {\mathbb Z}_{m},\qquad\alpha\in\left\{1,\ldots,m\right\}$$

\noindent and the basis of ${\cal H}^{\prime \prime }$ consist of:

$$\varphi^{\prime\prime}_{\beta}=e^{-\frac{1}{2}\sigma_{_{2}} y^{2}-c_{_{2}}x}\delta_{\beta}^{\nu }$$

\noindent where

$$\sigma_{_{2}} =i\tau_{2} \frac{l}{k-l\,\theta },\qquad \nu\in {\mathbb Z}_{l},\qquad\alpha\in\left\{1,\ldots,l\right\}$$
\noindent We assume that  ${\cal E}^{\prime }$, considered as a right $T_{\theta }$--module, is isomorphic to $E_{n,m}$ , and that ${\cal E}^{\prime\prime }$, considered as a left $T_{\theta }$--module, is isomorphic to $E^{\prime}_{k,l}$. We can use (\ref{prodbasis}) to calculate $\varphi^{\prime}_{\alpha}\underset{T_{\theta }}{\otimes }\varphi^{\prime\prime}_{\beta}$. The condition that complex structures on ${\cal E}^{\prime }$ and ${\cal E}^{\prime\prime }$ correspond to the same complex structure on ${\cal E}$ implies $\tau_{1}=\tau_{2}$. Using this, we can check that the theta--functions that appear in (\ref{prodbasis}) do not depend on $z$ in our case.  Applying (\ref{prodbasis}), we obtain that, $\varphi^{\prime}_{\alpha}\underset{T_{\theta }}{\otimes }\varphi^{\prime\prime}_{\beta}$ maps to $\Xi_{\alpha\,\beta}(z,\Delta)\in{\cal H}$ which is of the form:

$$\Xi_{\alpha\,\beta}(z,\Delta) = \sum_{\gamma} c_{\alpha\,\beta}^{\gamma}\,\varphi_{\gamma}(z,\Delta)$$

\noindent where

\begin{eqnarray*}
\varphi_{\gamma}(z,\Delta)&=&\exp\left[\frac{\pi\,i\,\tau(n\,l+m\,k)(n+m\,\theta)}{k-l\,\theta}z^{2}-(c_{_{1}}+c_{_{2}})(n+m\,\theta)z\right]\,\delta_{\gamma}^{\Delta}
\end{eqnarray*}

\noindent constitute a basis of ${\cal H}$ and the constants $c_{\alpha\,\beta}^{\gamma}$ are given by:

\begin{align}c_{\alpha\,\beta}^{\gamma}=\Theta(s,t)\cdot e^{K}\nonumber\end{align}

\noindent where $\Theta(s,t)$ a classical theta--function,

$$K=-\frac{\pi\,i\,\tau\,l}{m(n\,l+m\,k)}\,\gamma^{2}-\frac{\left(l(n+m\,\theta)c_{_{1}}-m(k-l\,\theta)c_{_{2}}\right)}{m(n\,l+m\,k)}\,\gamma+\pi\,i\,s\,q_{o}^{2}+2\,\pi\,i\,t\,q_{o}$$

\noindent and:

$$s=-m\,l\,(n\,l+m\,k)$$

$$t=\frac{l\,\tau\,\Delta-q_{o}(n\,l+m\,k)}{r} + \frac{c_{_{1}}\,l\,(n+m\,\theta)-c_{_{2}}\,m\,(k-l\,\theta)}{2\pi\,i\,r}$$

\subsection*{Acknowledgments}
We would like to thank Yu. Manin for his interest in our work, and for information that our results could be useful in number theory (Yu. Manin has explained the relation of our results to number theory in his recent preprint \cite{Ma4}).
One of authors (A. Sch.) appreciates the hospitality of IHES and ESI where
a part of this work was done.

\section{Appendix}

Our candidate function $h(z, \Delta)$ for the map (\ref{map}) is
given by :

\begin{eqnarray*}h(z,\Delta)= \sum_{q\in\mathbb Z}\left[f((n+m\,\theta)z-\frac{n+m\,\theta}{m}q+\frac{l(n+m\,\theta)}{m(n\,l+m\,k)}\Delta, -q+a\,\Delta)\right.\\ \left.\cdot g((n+m\,\theta)z+\frac{k-l\,\theta}{l}q-\frac{k-l\,\theta}{n\,l+m\,k}\Delta, q)\right]
\end{eqnarray*}
\noindent For reference, let us recall all the actions we will need. If $f(x,\mu)\in E_{n,m}(\theta)$, and $g(y,\nu)\in E^{\prime}_{k,l}(\theta)$, then:

$$U_{1}\,f(x, \mu)=f(x-\frac{n+m\,\theta}{m}, \mu-1)$$
$$ U_{2}\,f(x, \mu)=e^{2\,\pi\,i(x-\mu\,n/m)}\,f(x, \mu)$$
$$U_{1}\,g(y, \nu)=g(y-\frac{k-l\,\theta}{l}, \nu-1)$$
$$ U_{2}\,g(y, \nu)=e^{2\,\pi\,i(y-\nu\,k/l)}\,g(y, \nu)$$
$$Z_{1}\,f(x, \mu)=f(x-\frac{1}{m}, \mu-a)$$
$$ Z_{2}\,f(x, \mu)=e^{2\,\pi\,i(x/(n+m\,\theta)-\mu/m)}\,f(x, \mu)$$

\noindent Let us start by checking the identifications. Replacing
$f$ with $U_{1}\,f$ in the definition of $h$ we obtain:
\begin{align*} \sum_{q\in\mathbb Z}&\left[f((n+m\,\theta)z-\frac{n+m\,\theta}{m}q+\frac{l(n+m\,\theta)}{m(n\,l+m\,k)}\Delta-\frac{n+m\,\theta}{m}, -q+a\,\Delta-1)\right.\\ &\qquad\qquad \left.\cdot g((n+m\,\theta)z+\frac{k-l\,\theta}{l}q-\frac{k-l\,\theta}{n\,l+m\,k}\Delta, q)\right]
\intertext{Letting $q\to q-1$ gives:} &= \sum_{q\in\mathbb
Z}\left[f((n+m\,\theta)z-\frac{n+m\,\theta}{m}q+\frac{l(n+m\,\theta)}{m(n\,l+m\,k)}\Delta,
-q+a\,\Delta)\right.\\ & \qquad\qquad\left.\cdot
g((n+m\,\theta)z+\frac{k-l\,\theta}{l}q-\frac{k-l\,\theta}{n\,l+m\,k}\Delta-\frac{k-l\,\theta}{l},
q-1)\right]
\end{align*}

\noindent This means that the map (\ref{map}) is compatible with
the identification $$\left(U_{1}\,f\right)\otimes g \sim f \otimes
\left(U_{1}\,g\right)$$ Similarly, replacing $f$ with $U_{2}\,f$
in (\ref{map}) we get:

\begin{align*}
\sum_{q\in\mathbb Z}&\exp\left[2\,\pi\,i\left((n+m\,\theta)z-\frac{n+m\,\theta}{m}q+\frac{l(n+m\,\theta)}{m(n\,l+m\,k)}\Delta-\frac{n(-q+a\,\Delta)}{m}\right)\right]\,f\cdot g\\
&=\sum_{q\in\mathbb Z}\exp\left[2\,\pi\,i\left((n+m\,\theta)z-\theta\,q+\frac{l(n+m\,\theta)}{m(n\,l+m\,k)}\Delta-\frac{a\,n\,\Delta}{m}\right)\right]\,f\cdot g\\
\intertext{using the relation $a\,n=1+b\,m$ gives:}
&=\sum_{q\in\mathbb
Z}\exp\left[2\,\pi\,i\left((n+m\,\theta)z+\frac{k-l\,\theta}{l}q+\frac{k-l\,\theta}{n\,l+m\,k}\Delta-\frac{k}{l}q\right)\right]\,f\cdot
g
\end{align*}
 This means that the map (\ref{map}) is compatible with
the identification $$\left(U_{2}\,f\right)\otimes g \sim f
\otimes\left(U_{2}\,g\right)$$

\noindent So the map  (\ref{map}) that was defined originally as a
map of the tensor product $E_{n,m}\otimes_{\mathbb
C}E^{\prime}_{k,l}$ descends to a map of
$E_{n,m}\otimes_{T_{\theta}}E^{\prime}_{k,l}$.

Next, we need to check that the new product module is in fact
$E^{\prime}_{a\,k+b\,l,n\,l+m\,k}(\theta^{\prime})$. That is, we
first need to show that the shift $\Delta\to\Delta+n\,l+m\,k$
leaves $h(z,\Delta)$ invariant; therefore we can consider $\Delta$
as an element of $\mathbb Z_{n\,l+m\,k}$:

\begin{align*}
h&(z,\Delta+n\,l+m\,k)=\\ &\sum_{q\in\mathbb Z}\left[f((n+m\,\theta)z-\frac{n+m\,\theta}{m}q+\frac{l(n+m\,\theta)}{m(n\,l+m\,k)}(\Delta+n\,l+m\,k), -q+a\,(\Delta+n\,l+m\,k))\right.\\ & \qquad\left.\cdot g((n+m\,\theta)z+\frac{k-l\,\theta}{l}q-\frac{k-l\,\theta}{n\,l+m\,k}(\Delta+n\,l+m\,k), q)\right]\\
\end{align*}

\noindent Let $q\to q+l$. Then the extra terms in the first argument of $f$ are:

$$-\frac{n+m\,\theta}{m}l+\frac{l(n+m\,\theta)}{m(n\,l+m\,k)}(n\,l+m\,k)=0$$

\noindent The extra terms in the second argument of $f$ are:

\begin{align*}-l+a\,n\,l+a\,m\,k &\equiv -l-a\,n\,l \pmod m \\
& \equiv l\,(a\,n-1) \pmod m \\
& \equiv l\,b\,m \pmod m \\
& \equiv 0 \pmod m \\
\end{align*}

\noindent Similarly, the extra terms in the arguments of $g$ are respectively:

$$(k-l\,\theta)-(k-l\,\theta)=0$$
and
$$l\equiv 0\pmod l$$

\noindent It follows that:

$$h(z,\Delta+n\,l+m\,k)=h(z,\Delta)$$

\noindent as required. Finally we need to check that the action of
the endomorphism algebra generators $Z_{1}, Z_{2}$ on
$h(z,\Delta)$ induced by the action of $Z_{1}$ and $Z_{2}$ on $f$
describe the standard module
$E^{\prime}_{a\,k+b\,l,n\,l+m\,k}(\theta^{\prime})$:

\begin{align*}
Z_{1}&\,h(z,\Delta)\\
&= \sum_{q\in\mathbb Z}\left[f((n+m\,\theta)z-\frac{n+m\,\theta}{m}q+\frac{l(n+m\,\theta)}{m(n\,l+m\,k)}\Delta-\frac{1}{m}, -q+a\,\Delta-a)\right.\\ &\qquad\qquad \left.\cdot g((n+m\,\theta)z+\frac{k-l\,\theta}{l}q-\frac{k-l\,\theta}{n\,l+m\,k}\Delta, q)\right] \\
&= \sum_{q\in\mathbb Z}\left[f((n+m\,\theta)\left(z+\alpha\right)-\frac{n+m\,\theta}{m}q+\frac{l(n+m\,\theta)}{m(n\,l+m\,k)}(\Delta-1), -q+a\,(\Delta-1))\right.\\ &\qquad\qquad \left.\cdot g((n+m\,\theta)(z+\beta)+\frac{k-l\,\theta}{l}q-\frac{k-l\,\theta}{n\,l+m\,k}(\Delta-1), q) \right]\\
\end{align*}

\noindent where

\begin{align*}
\alpha &=\frac{1}{n+m\,\theta}\left(\frac{l(n+m\,\theta)}{m(n\,l+m\,k)}-\frac{1}{m}\right)\\
&=-\frac{k-l\,\theta}{(n+m\,\theta)(n\,l+m\,k)}\\
&= \beta
\intertext{and can be rewritten in a simpler way using the relation $a\,n-b\,m=1$:}
&=-\frac{a\,k+b\,l}{n\,l+m\,k}+\frac{b+a\,\theta}{n+m\,\theta}\\
&=-\frac{a\,k+b\,l}{n\,l+m\,k}+\hat\theta
\end{align*}

\noindent Therefore:
\begin{align*}
Z_{1}&\,h(z,\Delta)=Z_{1}\,h(z-\frac{a\,k+b\,l}{n\,l+m\,k}+\hat\theta,\Delta-1)
\end{align*}

\noindent as required for $E^{\prime}_{a\,k+b\,l,n\,l+m\,k}(\theta^{\prime})$. Similarly:

\begin{align*}
Z_{2}\,h(z,\Delta)&\\
&=\sum_{q\in\mathbb Z}\exp\left[2\,\pi\,i\left(z-\frac{q}{m}+\frac{l}{m(n\,l+m\,k)}\Delta-\frac{1}{m}(-q+a\,\Delta) \right) \right]f\cdot g\\
&=\sum_{q\in\mathbb Z}\exp\left[2\,\pi\,i\left(z+\frac{l(1-a\,n)-a\,m\,k}{m(n\,l+m\,k)}\Delta \right) \right]f\cdot g\\
\intertext{substituting $-b\,m$ in for $1-a\,n$ gives:}
&=\sum_{q\in\mathbb Z}\exp\left[2\,\pi\,i\left(z-\frac{a\,k+b\,l}{n\,l+m\,k}\Delta \right)\right]f\cdot g\\
&=\exp\left[2\,\pi\,i\left(z-\frac{a\,k+b\,l}{n\,l+m\,k}\Delta \right) \right]\,h(z,\Delta)
\end{align*}
also as required.

Our considerations prove that the formula (\ref{map}) specifies a map from $E_{n,m}(\theta)\underset{T_{\theta }}{\otimes }E_{k,l}^{\prime }(\theta)$ into $E_{a,k+b\,l,n\,l+m\,k}^{\prime}(\theta^{\prime})$ and this map is compatible with the structure of $T_{\theta^{\prime}}$ right modules on these objects. It is easy to check that this map is an isomorphism (either directly or using the general result about tensor products of modules over noncommutative tori).

One can use also the technique suggested in \cite{Sc2}; the latter approach has the advantage that it gives a regular way to obtain (\ref{map}). Namely we can start with the description of the linear space that is dual to the tensor product $E_{n,m}(\theta)\underset{T_{\theta }}{\otimes }E_{k,l}^{\prime }(\theta)$. We notice that a continuous linear functional on the tensor product of $E_{n,m}(\theta)$ and $E_{k,l}^{\prime }(\theta)$ over $\mathbb C$ is specified by a generalized function (a distribution) $\varphi(x,\mu,y,\nu)$ where $x, y \in \mathbb R, \mu\in \mathbb Z_{m}, \nu\in \mathbb Z_{l}$; it transforms $f\otimes g$ into:

\begin{eqnarray*}
\sum_{\mu\in\mathbb Z_{m}}\sum_{\nu\in\mathbb Z_{l}}\int f(x,\mu)\,g(y,\nu)\,\varphi(x,\mu,y,\nu)\,d\,x\,d\,y
\end{eqnarray*}

\noindent This functional descends to the tensor product over $T_{\theta}$ if:

\begin{align*}\varphi(x,\mu,y,\nu)\cdot e^{2\,\pi\,i(x-\mu\,n/m)}=\varphi(x,\mu,y,\nu)\cdot e^{2\,\pi\,i(y-\nu\,k/l)}\\ \intertext{and}
\varphi(x+\frac{n}{m}+\theta,\mu+1,y ,\nu)=\varphi(x,\mu,y+\frac{k}{l}-\theta,\nu+1)\end{align*}

Solving this system of equations for $\varphi$, we arrive at a description of the space that is dual to the tensor product we are interested in. Using this description, it is easy to obtain (\ref{map}).


\begin{thebibliography}{10}

\bibitem{Co1}
A.~Connes.
\newblock ${C}^{*}$--alg\`ebres et g\'eom\'etrie diff\'erentielle.
\newblock {\em C. R. Acad. Sci. Paris}, 290:599--604, 1980.

\bibitem{Co2}
A.~Connes.
\newblock {\em Noncommutative Geometry}.
\newblock Academic Press, 1994.

\bibitem{Co3}
A.~Connes and M.~Rieffel.
\newblock {Y}ang--{M}ills for noncommutative two--tori.
\newblock {\em Contemporary Math.}, 66:237--266, 1987.

\bibitem{Co4}
A.~Connes, M.~Douglas, and A.~Schwarz.
\newblock Noncommutative geometry and matrix theory: Compactification on tori.
\newblock {\em J. High Energy Phys.}, (2), 1998.

\bibitem{Ko1}
A.~Konechny and A.~Schwarz.
\newblock {BPS} states on noncommutative tori and duality.
\newblock {\em Nucl. Phys.}, B(550):561--584, 1999.

\bibitem{Ko2}
A.~Konechny and A.~Schwarz.
\newblock Moduli spaces of maximally supersymmetric solutions on noncommutative
  tori and noncommutative orbifolds.
\newblock {\em JHEP}, 09:1--23, 2000.

\bibitem{Ko3}
A.~Konechny and A.~Schwarz.
\newblock Introduction to {M}(atrix) theory and noncommutative geometry.
\newblock {\em arXiv:hep-th/0012145}.

\bibitem{Ma1}
Yu. Manin.
\newblock {\em Quantized theta--functions}, volume 102.
\newblock Common trends in Mathematics and Quantum Field Theories, 1990.

\bibitem{Ma2}
Yu. Manin.
\newblock Mirror symmetry and quantization of abelian varieties.
\newblock {\em arXiv:math.AG/0005143}.

\bibitem{Ma3}
Yu. Manin.
\newblock Theta functions, quantum tori and {H}eisenberg groups.
\newblock {\em Lett. Math. Phys.}, 56, June 2001.

\bibitem{Ma4}
Yu. Manin.
\newblock Real multiplication and noncommutative geometry.
\newblock {\em arXiv:math.AG/0202109}.

\bibitem{Mu1}
D.~Mumford.
\newblock {\em Tata Lectures on Theta {I}}.
\newblock Birkh\"{a}user, 1983.

\bibitem{Mu2}
D.~Mumford.
\newblock {\em Tata Lectures on Theta {III} (with {M}. {N}ori and {P}.
  {N}orman)}.
\newblock Birkh\"{a}user, 1991.

\bibitem{Ri1}
M.~Rieffel.
\newblock Projective modules over higher--dimensional noncommutative tori.
\newblock {\em Can. J. Math.}, XL(2):257--338, 1988.

\bibitem{Ri2}
M.~Rieffel and A.~Schwarz.
\newblock Morita equivalence of multidimensional noncommutative tori.
\newblock {\em Intl. J. of Math.}, 10(2):289--299, 1999.

\bibitem{Sc1}
A.~Schwarz.
\newblock Gauge theories on noncommutative spaces.
\newblock {\em ICMP lecture, arXiv:hep--th/0011261}.

\bibitem{Sc2}
A.~Schwarz.
\newblock Morita equivalence and duality.
\newblock {\em Nucl. Phys.}, B(534):720--738, 1998.

\bibitem{Sc3}
A.~Schwarz.
\newblock Theta--functions on noncommutative tori.
\newblock {\em Lett. Math. Phys.}, 58(1):81--90, 2001.

\end{thebibliography}
\end{document}